\documentclass[a4paper, 11pt, final]{wsheet}

\usepackage{amsmath}
\usepackage{amssymb}
\usepackage{oldmath}
\usepackage{enumerate}
\usepackage[mathscr]{eucal}
\usepackage{graphics}
\usepackage[dvips]{graphicx}

\usepackage{color}

\title[Weaved partial STS's]{Operation of weaving partial Steiner triple systems}

\author{Ma{\l}gorzata Pra{\.z}mowska, Krzysztof Pra{\.z}mowski}

\created{18.02.2010}
\modified{NOT}
\info{A continuation of polygon.tex, kulki.tex and muldipap.tex, ideas concerning iterating of drawing "centers" of cubes}

\def\konftyp(#1,#2,#3,#4){\left( {#1}_{#2}\, {#3}_{#4} \right)}

\def\pointnumbsymb{\mbox{\boldmath$\nu$}}
\def\pointnumbsymbmaly{\mbox{\scriptsize\boldmath$\nu$}}
\def\pointnumb{\mathchoice{\pointnumbsymb}{\pointnumbsymb}{\pointnumbsymbmaly}{\pointnumbsymbmaly}}
\def\linenumbsymb{\mbox{\boldmath$b$}}
\def\linenumbsymbmaly{\mbox{\scriptsize\boldmath$b$}}
\def\linenumb{\mathchoice{\linenumbsymb}{\linenumbsymb}{\linenumbsymbmaly}{\linenumbsymbmaly}}
\def\ranksymb{\mbox{\boldmath$r$}}
\def\ranksymbmaly{\mbox{\scriptsize\boldmath$r$}}
\def\rank{\mathchoice{\ranksymb}{\ranksymb}{\ranksymbmaly}{\ranksymbmaly}}
\def\granksymb{\mbox{\boldmath$\gamma$}}
\def\granksymbmaly{\mbox{\scriptsize\boldmath$\gamma$}}
\def\grank{\mathchoice{\granksymb}{\granksymb}{\granksymbmaly}{\granksymbmaly}}

\def\LineOn(#1,#2){\overline{{#1},{#2}\rule{0em}{1,5ex}}}
\def\triop{\circledcirc}

\def\GrasSpace(#1,#2){{\bf G}_{{#2}}({#1})}
\def\MuPap(#1){{\mathscr P}_{{#1}}}
\def\MuPapx(#1){{\mathscr M}_{{#1}}}
\def\mupapg(#1){G{\mathscr P}({#1})}

\def\CyclSpacex(#1,#2,#3){{#1}\circledast^{{#2}}{#3}}

\def\PolPap(#1,#2){\CyclSpacex({},#1,{#2})}

\def\lines{{\cal L}}
\def\razy{\bowtie}

\newenvironment{ctext}{%
  \par
  \smallskip
  \centering
}{%
 \par
 \smallskip
 \csname @endpetrue\endcsname
}

\def\myend{{}\hfill{\small$\bigcirc$}}

\newtheorem{reprx}[thm]{Representation}

\begin{document}

\maketitle

\begin{abstract}
  We introduce an operation of a kind of product which associates with a partial
  Steiner triple system another partial Steiner triple system, the starting one
  being a quotient of the result.
  We discuss relations of our product to some other product-like constructions
  and prove some preservation/non-preservation theorems.
  In particular, we show series of anti-Pasch Steiner triple systems
  which are obtained that way.
\strut\par\noindent{\it Key words}: 
  convolution (of a partial Steiner triple system and a group),
  Veblen (= Pasch) configuration,
  (partial) Steiner triple system,
  Desargues configuration,
  Fano configuration,
  (finite) affine space,
  affine slit space.
\end{abstract}

\section*{Introduction}

In the paper we introduce the operation of weaving, 
which associates with a partial Steiner triple system 
(shortly: with a PSTS)
$\goth M$ a ``product" $\PolPap(m,{\goth M})$  
of $\goth M$ and a cyclic group of order $m$ in such a way that 
a quotient of $\PolPap(m,{\goth M})$ wrt. to a congruence $\approx$ is $\goth M$,
and the coimage of a line of $\goth M$ under natural projection 
$\PolPap(m,{\goth M})\longrightarrow \PolPap(m,{\goth M})\diagup\approx$
is a generalization of the Pappus configuration.
Clarly, these properties do not characterize the operation of weaving uniquely,
and several constructions which have these properties can be found in the literature,
just to mention the operation of convolution and 
the ``product'' defined in \cite{griggs}, of an STS with a parallel class distinguished
and an abelian group.

\par
The notion of convolution was introduced in \cite{convol}
(in a slightly less general way than this adopted in this paper), 
though it was used, implicitly,
e.g. in \cite{hilbert}, \cite{antimitra}, \cite{antypasz1} (comp. also \cite{sparse}).
Both in the construction of the convolution and 
the construction of $m$-th weaved configuration
applied to a partial Steiner triple system $\goth M$
the constructed points are ``weighted" points of $\goth M$ i.e. pairs of the form
$a_i = (a,i)$, where $a$ is a point of $\goth M$ and $i$ is an element of 
a fixed (in case of weaving -- cyclic) group $\sf G$.
The lines are sets of triples of weighted points on lines of $\goth M$ whose
weights satisfy certain conditions. One can note an analogy between these constructions
and the product construction.
An analogy only, since the triples of weights on lines of the constructed configuration
need not yield any PSTS defined on $\sf G$.
%
%
Examples of classical configurations that are convolutions were already quoted in
\cite{convol}, these are e.g. the Veblen configuration 
(also called {\em the Pasch configuration}), 
the Reye configuration, the Pappus configuration.
Some of them are also weaved configurations. 
Some of them are also members of another family, family of configurations
presentable as cyclically inscribed triangles%
\footnote{%
This family contains all (i.e. three) $9_3$-configurations.%
}
(comp. \cite{polygon}).
In this paper we do not study this family on its own 
%
%
but its members play an important role in characterizations of weaved configurations
(cf. Lem. \ref{prop:bezpap} and its consequences).
%
%
The Pappus configuration is 
an important example of a configuration which is in each of these three families.
\par
We start with establishing general properties of weaved configurations.
The operation of weaving destroys most of classical configurations based on 
the Veblen configuration; in particular it gives 
Desargues-free (Prop. \ref{prop:des}) and 
Fano-free (Prop. \ref{prop:fana}) configurations.
They are also miter-free (Prop. \ref{prop:bezmitra}).
The operation of weaving applied to a Pasch-free configuration yields a 
Pasch-free configuration (Cor. \ref{cor:antypasz}). 
In particular, when applied to a Pasch-free
STS it produces a Pasch-free PSTS which has a unique, Pasch-free, completion
to an STS (Rem. \ref{rem:antypasz}).  
\par
Still, geometry of an $m$-th weaved configuration is relatively easy to understand.
In particular, it is easy to determine the triangles and cliques in it
(Lem. \ref{lem:triangles}), 
to characterize some characteristic subconfigurations (being direct generalizations
of the Pappus configuration, Prop. \ref{prop:pseudopapy}, the proof of \ref{prop:bezpap}), 
and to characterize (for $m > 3$) its automorphism group (Thm. \ref{thm:autpol}).
\par
After proving general properties of weaving operation, in the last section
we show some applications of the obtained results.
%
In particular, we get a method to obtain a class of STS's with parameters of
an affine space over $GF(3)$, 
which are not embeddable into any affine space over $GF(3)$
and which are Pasch-free (Prop. \ref{prop:afNOTaf}). 
Each line of the resulting configuration can be extended to an affine (sub)plane,
and planes are their maximal subspaces which are affine.
%

\section{Definitions and representations}\label{sec:defy}

Let ${\goth M} = \struct{S,\lines}$ be a partial Steiner triple system.
Let $m > 2$ be an integer.
Write 
  $X := S \times C_m$ 
and 
  ${\cal C} = \{ (i,j,k)\in C_m^3\colon i=j=k-1 \Lor i=k=j-1 \Lor j=k=i-1 \}$.
Finally, we define 
\begin{eqnarray}
  \lines_\ast & := & \left\{ \{(a,i),(b,j),(c,k)\} \colon
  \{ a,b,c \}\in\lines, (i,j,k)\in{\cal C} \right\},
  \\
  \PolPap(m,{\goth M}) & := & \struct{X,\lines_\ast}.
\end{eqnarray}
The structure $\PolPap(m,{\goth M})$ will be called a 
{\em configuration weaved from $\goth M$}
(more precisely, {\em the configuration $m$-weaved from $\goth M$}).
\par
Recall a similar construction of the {\em convolution} 
  ${\goth M}\razy_\varepsilon{\sf G}$
of an abelian group 
${\sf G} = \struct{G,0,+}$ and a partial Steiner triple system $\goth M$ (cf. \cite{convol}).
Let 
  $X := S \times G$,
  $\varepsilon\in G$, 
and
  ${\cal G}_{\varepsilon} = 
    \{(\alpha,\beta,\gamma)\in G^3 \colon \alpha+\beta+\gamma=\varepsilon \}$.
We set 
\begin{eqnarray}
  \lines_\varepsilon & := &
  \left\{ \{ (x,\alpha),(y,\beta),(z,\gamma) \} \colon \{x,y,z \}\in\lines,\,
  (\alpha,\beta,\gamma)\in{\cal G}_\varepsilon\right\},
  \\
  {\goth M}\razy_\varepsilon G & := & \struct{X, \lines_\varepsilon}.
\end{eqnarray}
\begin{fact}[{comp. \cite{twistfan}}]  
  Let $\varepsilon,e \in G$ and $f\in\Aut({\sf G})$. 
  Then
    ${{\goth M}\razy_{\varepsilon}{\sf G}} \cong {{\goth M}\razy_{\varepsilon + 3e}{\sf G}}
      \cong {{\goth M}\razy_{f(\varepsilon)}{\sf G}}$. 
  In particular,
    ${\goth M}\razy_{0} C_2 \cong {\goth M}\razy_{1} C_2$.
\end{fact}
The choice of the value `1' in the definition of the lines in $\lines_\ast$
may seem arbitrary, and one can consider the class
  $\big\{ \{ (a,i),(b,i),(c,i+\varepsilon) \} 
    \colon \{ a,b,c \}\in\lines,\; i\in C_m  \big\}$
with fixed $\varepsilon$ instead.
However, the obtained configuration may stay disconnnected,
and its connected components are isomorphic to 
$\PolPap(k,{\goth M})$ where $k$ is the rank of $\varepsilon$ in $C_m$.
Geometry of the convolution ${\goth M}\razy_\varepsilon{\sf G}$
may depend on $\varepsilon$, though.
\par
The structures of the form ${\goth M}\razy_0{\sf G}$ were studied in \cite{convol}
in much detail. We write, shortly, 
  ${\goth M}\razy{\sf G}$ instead of ${\goth M}\razy_0{\sf G}$.
Basic parameters of the structures defined above are easy to compute. 
\begin{fact}\label{fct:params}
  Let $\goth M$ be a $\konftyp(\pointnumb,\rank,\linenumb,3)$-configuration.
  Then $\PolPap(m,{\goth M})$ is a 
  $\konftyp(m \pointnumb,3\rank,3m\linenumb,3)$-con\-fi\-gu\-ra\-tion and
  ${\goth M}\razy_\varepsilon{\sf G}$ is a 
  $\konftyp( {\grank}\pointnumb,{\grank}{\rank}, {\grank}^2\linenumb, 3)$-con\-fi\-gu\-ra\-tion,
  where $\grank$ is the rank of $\sf G$.
  In particular, $\PolPap(m,{\goth M})$ and ${\goth M}\razy_\varepsilon{\sf G}$
  both are partial Steiner triple systems.
\end{fact}

Let us begin with some evident examples.
\begin{exm}
  Let $\goth T$ be a single-line structure i.e. let $|S| = 3$ and
    ${\goth T} = \struct{S,\{S\}}$.
  Then $\PolPap(m,{\goth T})$  
  is a series of cyclically inscribed triangles,
  as considered in \cite{polygon}. 
  In particular, 
  $\PolPap(3,{\goth T})$ is the Pappus configuration.
  It is known (cf. \cite{convol}) that 
  ${\goth T}\razy_0 C_3$ is the Pappus configuration as well.
\myend
\end{exm}
\begin{rem}
  For each partial Steiner triple system $\goth M$ we have
\begin{equation}\label{eq:polASconv}
   \PolPap(3,{\goth M}) \cong {\goth M}\razy_2 C_3 \cong {\goth M}\razy_1 C_3 
\end{equation}
\end{rem}

\begin{rem}\label{rem:commut}
  Let $m,n\geq 3$ be integers and 
  $\varepsilon$ be an element of an abelian group $\sf G$.
  Then 
    $\PolPap(m,{({\goth M}\razy_\varepsilon{\sf G})})
      \cong
      {\PolPap(m,{\goth M})}\razy_\varepsilon{\sf G}$
  and
  $
    \PolPap(m,{(\PolPap(n,{\goth M}))}) \cong
      \PolPap(n,{(\PolPap(m,{\goth M}))})
  $
  for each partial Steiner triple system $\goth M$.
\end{rem}
From \ref{fct:params} we get easily
\begin{exm}
  {\em There is no configuration $\goth M$
  and no element $\varepsilon$ of any group $\sf G$, 
  such that 
    $\PolPap(5,{\goth T}) \cong {\goth M}\razy_\varepsilon {\sf G}$}. 
\myend
\end{exm}
\begin{exm}
  $\PolPap(5,{\goth T})\not\cong {\goth T}\razy C_3$ and    
  $({\PolPap(5,{\goth T})})\razy C_3 \cong \PolPap(5,{({\goth T}\razy C_3)})$,
   so there is a structure which can be presented both
  as a convolution and as a weaved configuration
  of nonisomorphic structures.
\myend
\end{exm}

Finally, let us recall a few definitions from the general theory of
(partial) Steiner triple systems.
With each partial Steiner triple system ${\goth M} = \struct{S,\lines}$
we associate the partial binary operation $\triop$ defined on $S$ by the conditions:
\begin{eqnarray}
  p \triop q & := & \left\{
  \begin{array}{ll}
    p & \text{when } p = q
    \\
    r & \text{when } \{ p,q,r \}\in\lines.
  \end{array}
  \right. 
\end{eqnarray}
Let $\Delta = \{ p,q,r \}$ be a (nondegenerate) triangle in $\goth M$,
i.e. let it be a triple of pair wise collinear points not on a line.
We set 
\begin{equation}\label{def:Deltaprim}
  \Delta'  \; := \; \{ p\triop q, q\triop r, r \triop p \},
  \quad
  \Delta^{(n+1)} \; := \; {\Delta^{(n)}}',
\end{equation}
for each integer $n$ such that the points in $\Delta^{(n)}$ are pairwise collinear.
The structure $\goth M$ is called {\em Moufangian} iff 
$\Delta'$ is a line of $\goth M$ for every triangle $\Delta$.
The algebraic counterpart of this property expressed 
in terms of the partial algebra $\struct{S,\triop}$ 
is read as the known Moufang axiom:
\begin{ctext}\refstepcounter{equation}\label{tozs:muf}$
\strut\hfill  (p\triop q)\triop(p\triop r) = q\triop r\hfill\eqref{tozs:muf}
$ 
\end{ctext}
valid for every triangle $\{ p,q,r \}$.
Note that the Moufang property implies the Veblen axiom.

\section{General, subconfigurations}

In most parts proofs of evident statements are omitted.

Let us fix a partial Steiner triple system ${\goth M} = \struct{S,\lines}$.

\begin{lem}\label{lem:colin}
  Distinct points $(a,i),(b,j)$ of $\PolPap(m,{\goth M})$
  are collinear iff $a,b$ are distinct and collinear in $\goth M$ and
  $j \in \{ i,i+1,i-1 \}$.
\end{lem}

\begin{lem}\label{lem:triangles}
  Let $\Delta$ be a triangle in ${\goth B} = \PolPap(m,{\goth M})$.
  Then one of the following holds.
  \begin{sentences}\itemsep-2pt
  \item\label{triangles:typ1}
    $\Delta = \{ (a,i), (b,i), (c,i) \}$ for $i\in C_m$ 
    and a triangle $\Delta_0 := \{ a,b,c \}$ of $\goth M$.
    The set $\Delta'$ is a triangle in $\goth B$ if 
    $\Delta'_0$ is either a triangle or a line of $\goth M$ 
    (i.e. $\Delta_0$ is Moufangian).
    In the second case the (periodic) series
    $\Delta',\Delta'',\ldots,\Delta^{(m)},\Delta^{(m+1)}=\Delta'$
    consists of triangles distinct from $\Delta$.
    In any case $\Delta'$ is not a line of $\goth B$.
   \par
    Let $\Delta$ be a triangle of the form \eqref{triangles:typ1} and let $\Delta_0$ be
    a triangle as well. 
    Then $\Delta,\Delta',... ,\Delta^{(m-1)},\Delta^{(m)} = \Delta$
    is a periodic series iff there is an integer $m_0$ such that 
    $\Delta_0,\Delta'_0,...,\Delta^{(m_0-1)}_0,\Delta^{(m_0)}_0=\Delta_0$
    is a periodic series of triangles in $\goth M$ and $m_0$ divides $m$.
  \item\label{triangles:typ2}
    $\Delta = \{ (a,i), (b,i), (c,i) \}$ for $i\in C_m$ 
    and a line $\{ a,b,c \}$ of $\goth M$.
    Then the (periodic) series
    $\Delta,\Delta',\ldots,\Delta^{(m-1)},\Delta^{(m)} = \Delta$
    consists of triangles.
  \item\label{triangles:typ31}
    $\Delta = \{ (a,i), (b,i), (c,i-1) \}$ for $i\in C_m$, 
    where $\{ a,b,c \}$ is a line of $\goth M$.
    Then $\Delta'$ is not a line. For $m\neq 3$ it is not  a triangle
    (it consists of a pair of collinear points $p,q$ and a point $r$
    not collinear with any of $p,q$).
    For $m = 3$, $\Delta'$ is a triangle.
  \item\label{triangles:typ3}
    $\Delta = \{ (a,i), (b,i), (c,i-1) \}$ for $i\in C_m$, 
    where $\Delta_0 = \{ a,b,c \}$ is a triangle of $\goth M$.
    Then $\Delta'$ is not a line. If $m\neq 3$ then $\Delta'$ is not a triangle:
    it has either the form of \eqref{triangles:typ31} or
    consists of a triple of pairwise noncollinear points.
    If $m = 3$ then $\Delta'$ is a triangle iff the points in $\Delta'_0$
    are pairwise collinear.
  \item\label{triangles:typ4}
    $\Delta = \{ (a,i), (b,i), (c,i+1) \}$ for $i\in C_m$ 
    and a triangle $\Delta_0:=\{ a,b,c \}$ of $\goth M$.
    Then $\Delta'$ is a line of $\goth B$ iff $\Delta'_0$ is a line of 
    $\goth M$ (i.e. $\Delta_0$ yields a Veblen figure in $\goth M$).
    If $\Delta'_0$ is a triangle then $\Delta'$ is a triangle as well.
  \item\label{triangles:typ1:3}
    $\Delta = \{ (a,i), (b,i+1), (c,i+2) \}$ for $i\in C_m$
    and a triangle $\{ a,b,c \}$ of $\goth M$. The sequence $\Delta^{(j)}$
    is as in \eqref{triangles:typ1}.
  \item\label{triangles:typ2:3}
    $\Delta = \{ (a,i), (b,i+1), (c,i+2) \}$ for $i\in C_m$
    and a line $\{ a,b,c \}$ of $\goth M$. The sequence $\Delta^{(j)}$
    is as in \eqref{triangles:typ2}.
  \end{sentences}
  If $m\neq 3$ then $\goth B$ does not contain
  triangles of type \eqref{triangles:typ1:3} and
  \eqref{triangles:typ2:3}.
  If $m =3$ then triangles of these types may occur.
  If $\Delta$ has type \eqref{triangles:typ2:3} then there is in $\goth B$ a triangle
  $\Delta_1$ of type \eqref{triangles:typ2} such that 
    $\Delta\cup\Delta'\cup\Delta'' = \Delta_1\cup\Delta'_1\cup\Delta''_1$.
\end{lem}
\begin{rem}\small\normalfont
  Assume that $\goth M$ contains a triangle $\Delta_0$ such that 
    $\Delta_0,\Delta'_0,...,\Delta^{(m_0-1)}_0,\Delta^{(m_0)}_0=\Delta_0$
  is a periodic series of distinct triangles for some integer $m_0$.
  Let $n = LCM(m_0,m)$.
  Then $\PolPap(m,{\goth M})$ contains a triangle $\Delta$ such that
  a series of distinct triangles 
    $\Delta,\Delta',...,\Delta^{(n-1)},\Delta^{(n)}$ 
  exists and 
    $\Delta^{(n)} = \Delta$.
  Indeed, apply the construction of \eqref{triangles:typ1}.
  The construction of \eqref{triangles:typ4} yields a series with $n= m_0$.
\end{rem}

The only Veblen subconfigurations of $\PolPap(m,{\goth M})$
are determined by triangles characterized in \ref{lem:triangles}\eqref{triangles:typ4}.
As a direct consequence we get
\begin{cor}\label{cor:vebleny}
  A Veblen configuration contained in $\PolPap(m,{\goth M})$
  has form 
  \begin{ctext}
    $(a,i),(b,i),(c,i+1),(a\triop c,i),(b\triop c,i),(a\triop b,i+1)$,
  \end{ctext}
  where a triangle $\{ a,b,c \}$ yields a Veblen subconfiguration in $\goth M$ 
  and $i\in C_m$.
\end{cor}
\begin{cor}\label{cor:antypasz}
  If $\goth M$ does not contain any Veblen subconfiguration (is {\em Pasch-free}) then
  $\PolPap(m,{\goth M})$ is Pasch-free as well.
\end{cor}

An important consequence of \ref{cor:vebleny} is
\begin{prop}\label{prop:fana}
  The structure $\PolPap(m,{\goth M})$ does not contain any Fano subconfiguration.
  Actually, it is {\em anti-Fano}: no three diagonal points of a quadrangle
  contained in $\PolPap(m,{\goth M})$ are on a line.
\end{prop}
\begin{proof}
  Suppose that $\Delta$ is a triangle in $\PolPap(m,{\goth M})$ which spans
  a Fano configuration. In particular, $\Delta$ yields a Veblen configuration
  so, in view of \ref{cor:vebleny}, 
    $\Delta = \{(a,i),(b,i),(c,i+1)\}$
   where $\{ a,b,c \}$ is a triangle in $\goth M$ and $i\in C_m$.
   Write $b' = a\triop c$, $a' = b\triop c$, and $c' = a\triop b$,
   Then a quadrangle which spans a Fano plane has form
     $Q = \{ (a,i),(b,i),(a',i),(b',i) \}$, 
   provided the points $a',b',c'$ are collinear in $\goth M$. 
   Two of the diagonal points of $Q$ are $(c,i+1)$ and $(c',i+1)$.
   To get the third diagonal point of $Q$ we need a common point $p$
   of the lines through $a,a'$ and through $b,b'$. But then this third 
   diagonal point of $Q$ is $(p,i+1)$, and the points 
   $(c,i+1)$, $(c',i+1)$, and $(p,i+1)$ are never on a line of $\PolPap(m,{\goth M})$.
\end{proof}
\begin{rem}[ad the proof of \ref{prop:fana}]\label{rem:fana}
  Note that if the quadrangle $a,b,a',b'$ yields a Fano plane then 
  the diagonal points of the quadrangle $Q$ yield a nondegenerate
  triangle in $\PolPap(m,{\goth M})$.
\end{rem}
\begin{cor}\label{cor:emb2fano}
  Assume that $\goth M$ contains a Fano subconfiguration.
  Then $\PolPap(m,{\goth M})$ cannot be embedded to any projective space
  $PG(n,p)$ with even $p$ and $n \geq 2$.
\end{cor}
In an analogous fashion we get
\begin{prop}\label{prop:des}
  The structure $\PolPap(m,{\goth M})$ does not contain any Desargues subconfiguration.
  Actually, it is {\em anti-Desarguesian}: no three focuses of two perspective triangles
  contained in $\PolPap(m,{\goth M})$ are collinear.
\end{prop}
\begin{proof}
  In view of \ref{lem:triangles} two triangles which have a perspective center
  such that their focuses exist (i.e. corresponding sides of the triangles intersect
  in pairs) have form 
    ${\sf T}_1 := \left( (b,i),(c,i),(d,i) \right)$ and
    ${\sf T}_2 := \left( (a\triop b,i),(a\triop c,i),(a\triop d,i) \right)$,
  where 
    $(b, c, d)$ and
    $(a\triop b, a\triop c, a\triop d)$
  is a pair of triangles of $\goth M$ with a perspective center $a$.
  Then $(a,i+1)$ is the perspective center of ${\sf T}_1$ and ${\sf T}_2$.
  The focuses of ${\sf T}_1$ and ${\sf T}_2$ are the points
    $(b\triop c,i+1)$, $(c\triop d,i+1)$, and $(d\triop b,i+1)$.
  These points are not collinear.
\end{proof}
\begin{rem}[ad the proof of \ref{prop:des}]\label{rem:des}
    Note that if the triangles 
      $(b, c, d)$ and
      $(a\triop b, a\triop c, a\triop d)$
    yield a Desargues configuration i.e. their three focuses colline
    then the focuses of ${\sf T}_1$ and ${\sf T}_2$ yield a nondegenerate
    triangle in $\PolPap(m,{\goth M})$.
\end{rem}

As a direct consequence of \ref{prop:des} and, in particular, \ref{rem:des},
we get
\begin{cor}\label{cor:emb2des}
  Assume that $\goth M$ contains a Desargues subconfiguration.
  Then $\PolPap(m,{\goth M})$ cannot be embedded to any Desarguesian projective space.
 \par
  Consequently, neither $\PolPap(m,{\GrasSpace(X,2)})$
  for any $X$ with $|X|\geq 5$ 
 \footnote{\label{foot:gras}%
  The points of the incidence structure $\GrasSpace(X,2)$ are the two-element
  subsets of a set $X$, and the lines are the three-element subsets of $X$,
  the incidence being the inclusion
  (cf. \cite{maldem}, \cite{perspect}). These structures constitute a class
  of copolar spaces (cf. \cite{hall:copol}) not associated with any quadric.%
  },
  nor $\PolPap(m,{PG(n,2)})$ for any $n\geq 3$
  can be embedded into a Desarguesian projective space for any $m\geq 3$.
\end{cor}

In essence, \ref{prop:des} can be generalized to a wider class of ${10}_{3}$-configurations.
\begin{rem}\label{rem:K4}
  Let $\goth K$ be a ${10}_{3}$-configuration that contains a Veblen subconfiguration.
  The structure $\PolPap(m,{\goth M})$ does not contain any subconfiguration
  isomorphic to $\goth K$ for any $m > 3$ and any partial Steiner triple 
  system $\goth M$.
\end{rem}
\begin{proof}
  In accordance with \cite{klik:VC} there are exactly $6$ configurations $\goth K$
  of the form considered in \ref{rem:K4} and each one can be presented as
  a ``closure" of a $K_4$-graph. That means $\goth K$ contains a $K_4$-graph
  and ``third points" on the edges of this graph yield a Veblen subconfiguration.
  Suppose that ${\goth B} := \PolPap(m,{\goth M})$ contains $\goth K$.
  Let $\goth V$ be the respective Veblen subconfiguration of $\goth K$ with the points
  (cf. \ref{cor:vebleny}) 
    $(p,i_0),(q,i_0),(r,i_0+1)$ (on a line) and $(a,i_0),(b,i_0),(c,i_0+1)$ (a triangle).
  The edges of $\goth G$ which pass through 
  $(c,i_0+1)$ have form 
    $(x,i_0),(y,i_0)$ (1) or $(x,i_0+1),(y,i_0+2)$ (2), and
  through $(r,i_0+1)$ have form
    $(z,i_0),(t,i_0)$ (3) or $(z,i_0+1),(t,i_0+2)$ (4) resp.
  Suppose (1) \& (3). Without loss of generality we can assume $x\neq t$ and then
    $(x,i_0)\triop(z,i_0) = (x\triop z, i_0+1)$
  is another point of $\goth V$, which is impossible.
  Other cases are considered analogously.
\end{proof}

Another configuration frequently considered in combinatorics of STS's
is the {\em miter configuration}
(cf. e.g. \cite{antimitra}, where anti-miter STS's were studied). 
A triangle $\{ a,b,c \}$
determines a miter configuration with the center $a$ when the equality
\begin{ctext}\refstepcounter{equation}\label{tozs:mitra}
\strut\hfill $a\triop(b \triop c) = (a\triop b)\triop(a\triop c)$
\hfill\eqref{tozs:mitra}
\end{ctext}
holds, and the configuration in question consists of the points
  $a, b, c, a\triop b, a\triop c, b \triop c, a\triop(b \triop c)$.
\begin{prop}\label{prop:bezmitra}
  The structure $\PolPap(m,{\goth M})$ does not contain any miter-configuration.
\end{prop}
\begin{proof}
  Analyzing all the possibilities given in \ref{lem:triangles} we check that
  no triangle in $\PolPap(m,{\goth M})$ may satisfy equation \eqref{tozs:mitra}.
\end{proof}
\begin{cor}\label{cor:afNOTemb}
  Neither the M{\"o}bius\  ${8}_{3}$-configuration (cf. \cite{hilbert})
  nor the affine plane $AG(2,3)$ 
  can be embedded into a weaved configuration $\PolPap(m,{\goth M})$.
\end{cor}
\begin{proof}
  It suffices to note that the ${8}_{3}$-configuration results from 
  $AG(2,3)$ by removing a point and all the lines through it,
  and the miter configuration results from the ${8}_{3}$-configuration
  by omitting a point and all the lines through it.
\end{proof}

In the sequel we need criterions
which enable us to distinguish triangles of 
form \ref{lem:triangles}\eqref{triangles:typ1}
and those of form \ref{lem:triangles}\eqref{triangles:typ2}.
To this aim we must recall a fragment of \cite{polygon}.
Let us start with a naive approach. Consider a triangle
$\Delta^{(0)}$. Inscribe a triangle $\Delta^{(1)}$ into $\Delta^{(0)}$.
Inductively, inscribe a triangle $\Delta^{(i+1)}$ into $\Delta^{(i)}$.
Continue this procedure $(m-1)$ times so as a triangle $\Delta^{(m-1)}$
is obtained. Finally, inscribe $\Delta^{(0)}$ into $\Delta^{(m-1)}$.
\par
The obtained configuration is uniquely determined by
a permutation $\gamma$ of $C_3$ so as (up to an isomorphism)
the points of the arising configuration denoted by $\Pi^m_\gamma$
are the elements of $C_3\times C_m$ and 
the lines of $\Pi^m_\gamma$ are the sets
  $\left\{ (a,i), (b,i), (c,i+1) \right\}$ for all the triples $a,b,c$ such that
  $C_3 = \{ a,b,c \}$ and $i=0,1,\ldots,m-2$, 
and the sets
  $\left\{ (a,m-1),(b,m-1),(\gamma(c),0) \right\}$ with $C_3 = \{ a,b,c \}$
(cf. \cite{polygon}, slightly modified).
For fixed integer $m$ there are up to an isomorphism exactly
three configurations of the form $\Pi^m_\gamma$:
with 
  $\gamma = \id$, $\gamma = \tau_1$, and $\gamma = \sigma_0$ 
  ($\tau_1$ is the translation on $1$: $\tau_1(a) = a+1$, and
   $\sigma_0$ is the reflection in $0$: $\sigma_0(a) = -a$).
Recall also that there are exactly three $\konftyp(9,3,9,3)$-configurations
and these are 
$\Pi^3_{\id}$ (= the Pappus configuration),
$\Pi^3_{\tau_1}$, and $\Pi^3_{\sigma_0}$.
The following simple observation shows 
a close connection between weaved configurations and series of inscribed 
triangles (cf. \ref{lem:triangles}\eqref{triangles:typ2}).
\begin{equation}
  \Pi^m_{\id} \cong \PolPap(m,{\goth T}) \text{ for each integer } m \geq 3.
\end{equation}

Contrary to \ref{prop:fana} and \ref{prop:des}, from \ref{lem:triangles}\eqref{triangles:typ4}
we have immediately
\begin{prop}\label{prop:pseudopapy}
  If $\goth M$ contains $\Pi^k_\gamma$ for some permutation $\gamma$ of $C_3$
  and some integer $k$ then $\PolPap(m,{\goth M})$ also contains $\Pi^k_\gamma$.
  In particular, if $\goth M$ contains a Pappus subconfiguration then 
  $\PolPap(m,{\goth M})$ also contains a Pappus subconfiguration
  (comp. \ref{cor:afNOTemb}).
\end{prop}
\begin{proof}
  Let $\Delta_0 = \{ a,b,c \}$ be a triangle of $\goth M$ which determines a cyclic
  series 
    $\Delta'_0,...,\Delta^{(k)}_0 = \Delta_0$ 
  of inscribed triangles 
  such that 
    ${\cal X} = \bigcup \{\Delta^{(j)}_0\colon j=0,...,k-1\}$ 
  yields the $\Pi^k_\gamma$ subconfiguration of $\goth M$.
  Take any $i\in C_m$ and set 
    $\Delta = \{ (a,i),(b,i),(c,i+1) \}$.
  It is seen that the triangle $\Delta$ yields in $\PolPap(m,{\goth M})$ a cyclic
  series of length $k$ of inscribed triangles such that the set
    $\bigcup \{\Delta^{(j)}\colon j=0,...,k-1\} \subset {\cal X}\times\{ i,i+1 \}$ 
  yields the $\Pi^k_\gamma$ subconfiguration of $\PolPap(m,{\goth M})$.
\end{proof}
With a bit more subtle analysis we can also prove
\begin{rem}\label{rem:papax}
  Assume that $\goth M$ satisfies the projective Pappus axiom.
  Then every three diagonal points of a hexagon of $\PolPap(m,{\goth M})$
  inscribed into two lines, are on a line.
\end{rem}
\begin{proof}
  Let $p_1,\ldots,p_6$ be a hexagon of $\PolPap(m,{\goth M})$ inscribed into two lines:
  i.e. assume that $\{ p_1,p_3,p_5 \}$ and $\{ p_2,p_4,p_6 \}$ are two lines.
  Let the corresponding diagonal points be:
    $q_1$ on $\LineOn(p_1,p_2),\LineOn(p_4,p_5)$,
    $q_2$ on $\LineOn(p_2,p_3),\LineOn(p_5,p_6)$, and
    $q_3$ on $\LineOn(p_3,p_4),\LineOn(p_6,p_1)$.
  Then $p_2,p_3,p_4$ is a triangle inscribed into the triangle $p_1,p_5,q_1$,
  and $q_2,p_6,q_3$ are third points on the sides of the triangle 
  $p_2,p_3,p_4$. Analyzing possible ways in which series of inscribed triangles
  may be obtained, with the help of \ref{lem:triangles} we note that the triples
  $p_5,p_6,q_2$ and $p_1,p_6,q_3$ are collinear only in case \eqref{triangles:typ4}
  and in that case 
    $p_1 = (a_1,i)$, $p_2= (a_2,i)$, $p_3 = (a_3,i+1)$, $p_4 = (a_4,i)$,
    $p_5 = (a_5,i)$, $p_6 = (a_6,i+1)$, 
  where $a_1,\ldots,a_6$ is a hexagon in $\goth M$ 
  with the diagonal points $b_1,b_2,b_3$ such that 
    $q_1 = (b_1,i+1)$, $q_2 = (b_2,i)$, $q_3 = (b_3,i)$. 
  Now the claim is evident.
\end{proof}

\section{Automorphisms}

\begin{lem}\label{lem:aut:evident}
  Let $u \in C_m$ and
  $f \in\Aut({\goth M})$.
  Then
    the map $f\times \tau_u \colon S\times C_m \ni(a,i)\longmapsto (f(a),i+u)$
    is an automorphism of $\PolPap(m,{\goth M})$.
\end{lem}

As a simple consequence of \ref{lem:triangles} we get
\begin{lem}\label{prop:bezpap}
  Assume that no one of the following configurations is contained in $\goth M$:
  \def\dziel{\mathrel{\vert}}
  \begin{enumerate}[\rm(i)]\itemsep-2pt
  \item
    $\Pi^{m_0}_{\id}$, where $m_0 \dziel m$,
  \item
    $\Pi^{m_0}_{\sigma_0}$, where $(2 m_0) \dziel m$, and
  \item
    $\Pi^{m_0}_{\tau_1}$, where $(3 m_0) \dziel m$.
  \end{enumerate}
  The family of sets $L\times C_m$ with $L$ ranging over the lines of 
  $\goth M$ is definable in terms of the geometry of $\PolPap(m,{\goth M})$.
\end{lem}
\begin{proof}
  We need to provide an analysis of triangles 
  slightly more subtle than in \ref{lem:triangles}.
  With a triangle $\delta = (a,b,c)$ (a sequence, not a set!) we associate 
  {\em the sequence} $\delta' = (a\triop b,b\triop c,c\triop a)$. 
  As in \eqref{def:Deltaprim}
  we introduce the symbols $\delta^{(i)}$.
  We claim that the following conditions are equivalent
  \def\dziel{\mathrel{\vert}}
  \begin{enumerate}[(a)]\itemsep-2pt
  \item\label{wrr1}
    ${\cal X} = L \times C_m$ for a line $L$ of $\goth M$,
  \item\label{wrr2}
    ${\cal X} = \Delta \cup \Delta' \cup\ldots\cup\Delta^{(m-1)}$ for a triangle
    $\Delta= \{ p,q,r \}$ of $\PolPap(m,{\goth M})$ and $\delta = (p,q,r)$ such that 
    $\Delta',\ldots,\Delta^{(m-1)}$ consists of distinct triangles,
    ${(\Delta^{(m-1)})}' = \Delta^{(m)} = \Delta$, and $\delta^{(m)} = \delta$.
  \end{enumerate}
  for any set $\cal X$ of points of $\PolPap(m,{\goth M})$.
  \par\noindent
  Implication \eqref{wrr1}$\Rightarrow$\eqref{wrr2} is a direct consequence
  of \ref{lem:triangles}\eqref{triangles:typ2}.
  Assume \eqref{wrr2} and suppose that \eqref{wrr1} is not valid;
  in view of \ref{lem:triangles} we get that one of the following holds:
  \ref{lem:triangles}\eqref{triangles:typ1},
  \ref{lem:triangles}\eqref{triangles:typ4}, or $m =3$ and
  \ref{lem:triangles}\eqref{triangles:typ1:3}, 
  \ref{lem:triangles}\eqref{triangles:typ3}, 
  \ref{lem:triangles}\eqref{triangles:typ31}, or
  \ref{lem:triangles}\eqref{triangles:typ2:3}.
  \par
  In the first four cases $\Delta$ is associated with a triangle $\Delta_0 = \{ a,b,c \}$
  of $\goth M$. Write $\delta_0 = (a,b,c)$. 
  It is seen that $\Delta^{(i+1)}$ can be considered as a triangle 
  inscribed into $\Delta^{(i)}$. 
  From \eqref{wrr2} we get that $\Delta_0^{(m)} = \Delta_0$.
  Let $m_0$ be the least integer with $\Delta_0^{(m_0)} = \Delta_0$.
 \par
  Suppose, first, that
    \ref{lem:triangles}\eqref{triangles:typ1} holds.
  Clearly, $m_0 \dziel m$ i.e. $m = m_0 k$ for some integer $k$. 
  There is a permutation $\gamma$ of $C_3$ such that
    $\delta_0^{(m_0)} = (\gamma(a),\gamma(b),\gamma(c))$
  and then the points of 
    $\bigcup \{ \Delta_0^{(i)} \colon i = 0,...,m_0-1 \}$
  yield in $\goth M$ the configuration $\Pi^{m_0}_\gamma$.
  From assumption, $\delta_0^{(m)} = \delta_0$ and thus $\gamma^{k} = \id$.
  If $\gamma$ is a translation then $3 \dziel k$ and if $\gamma$ is a 
  reflection then $2 \dziel k$, which contradicts assumptions.
 \par
  Next, suppose that 
   \ref{lem:triangles}\eqref{triangles:typ4} holds.
  In that case from \eqref{wrr2} we get $m_0 = m$, and then $\delta^{(m)} = \delta$
  yields that $\goth M$ contains $\Pi^{m_0}_{\id}$, which is impossible.
 \par
  Now, let $m = 3$. 
  Assume that
   \ref{lem:triangles}\eqref{triangles:typ3} or 
   \ref{lem:triangles}\eqref{triangles:typ1:3} holds. 
  Clearly, $m_0\dziel m$ and thus $m_0 = 3$.
  From $\delta^{(3)} = \delta$ we get that $\Pi^3_{\id}$ is contained in $\goth M$,
  which contradicts assumptions.
\par 
  Finally, consider the last two cases i.e. assume that
  \ref{lem:triangles}\eqref{triangles:typ31} or 
  \ref{lem:triangles}\eqref{triangles:typ2:3} holds. 
  In these cases $\Delta$ arises from a  line $L$ of $\goth M$.
  It is seen that $\bigcup_{j=0}^{2} \Delta^{(j)} = L\times C_3$ i.e. \eqref{wrr1} holds.
\end{proof}
A structure $\goth M$ which satisfies the assumptions of \ref{prop:bezpap}
will be called {\em anti-$m$-po\-ly\-pap\-pian}.
It is seen that a Moufangian configuration $\goth M$ 
is anti-$m$-polypappian for each integer $m \geq 3$.

\begin{cor}\label{cor:aut:iloraz}
  Let $\goth M$ be an anti-$m$-polypappian PSTS with point degree $>1$.
  To every 
    $F\in\Aut({\PolPap(m,{\goth M})})$
  there corresponds a map $\alpha_F \in\Aut({\goth M})$ such that 
    $F(L\times C_m) = \alpha_F(L)\times C_m$ 
  for every line $L$ of $\goth M$.
\end{cor}

From \ref{cor:aut:iloraz}, in view of \ref{lem:aut:evident}
to characterize the group 
$\Aut({\PolPap(m,{\goth M})})$ where 
$\goth M$ is anti-$m$-polypappian and has the points of degree $>1$
it suffices to determine the kernel of the epimorphism 
\begin{ctext}
  $\alpha\colon \Aut({\PolPap(m,{\goth M})}) \longrightarrow \Aut({\goth M})$.
\end{ctext}

\begin{lem}\label{lem:pom}
  Assume that $m>3$.
  Let $F\in\Aut({\PolPap(m,{\goth M})})$ 
  such that $F(\{ a \}\times C_m) = \{ a \} \times C_m$ for each point $a$ of $\goth M$.
  Let $L = \{ a,b,c \}$ be a line of $\goth M$ and 
  $F(a,i_0) = (a,j_0)$ for some $i_0,j_0\in C_m$.
  \begin{sentences}\itemsep-2pt
  \item\label{pom:war1}
    Then $F(b,i_0)=(b,j_0)$ and $F(c,i_0)=(c,j_0)$.
  \item\label{pom:war2}
    Moreover, $F(x,i) = (x,i+(j_0-i_0))$ for each $x\in\{ a,b,c \}$.
  \end{sentences}
\end{lem}
\begin{proof}
  Let us write $\Delta(L,i) = \{ (a,i),(b,i),(c,i) \}$ for $i\in C_m$ and
  a line $L = \{ a,b,c \}$ of $\goth M$.
 \par
  Clearly, $F$ preserves the class of triangles of 
  $\PolPap(m,{\goth M})$ of the form \eqref{triangles:typ2} of \ref{lem:triangles}
  and thus $F(\Delta(L,i_0)) = \Delta(L,j)$ for some $j$. 
  From assumptions, $j=j_0$.
  This justifies \eqref{pom:war1}.
 \par
  Note that $\Delta(L,i)' = \Delta(L,i+1)$ for each $i\in C_m$.
  Therefore, 
  \begin{math}
   F(\Delta(L,i_0+1)) = F(\Delta(L,i_0)') = (F(\Delta(L,i_0)))' =
   \Delta(L,j_0)' = \Delta(L,j_0+1),
  \end{math}
  which gives $F(a,i_0+1) = (a,j_0+1)$.
  Inductively, we get $F(a,i_0+v) = (a,j_0+v)$ for each $v \in C_m$.
  This proves \eqref{pom:war2}.
\end{proof}

\begin{lem} 
  Assume that $\goth M$ is anti-$m$-polypappian, connected, 
  and with the points of degree $>1$.
  Let $F\in\ker(\alpha)$ and $m > 3$.
  Then there is $u\in C_m$ such that $F = \id \times\tau_u$.
\end{lem}
\begin{proof}
  By \ref{lem:pom}\eqref{pom:war1}, 
  for each line $L$ of $\goth M$ there is a bijection $\beta_F^L$ of $C_m$
  such that 
    $F(a,i) = (a,\beta_F^L(i))$ for every point $a$ on $L$.
  From the connectedness, 
    $\beta_F^{L'} = \beta_F^{L''}$ for any two lines $L',L''$ of $\goth M$. 
  Thus there is a bijection $\beta_F$ of $C_m$ such that 
    $F(a,i) = (a,\beta_F(i))$ for every point $a$ of $\goth M$.
  From \ref{lem:pom}\eqref{pom:war2} we get that 
    $\beta_F = \tau_u$ for some $u\in C_m$.
\end{proof}

As an immediate corollary we get
\begin{thm}\label{thm:autpol}
  Let $m > 3$ and let $\goth M$ be an anti-$m$-polypappian
  connected partial Steiner triple system
  with the points of degree $>1$.
  Then 
  \begin{equation}
    \Aut({\PolPap(m,{\goth M})}) = 
    \left\{ f\times\tau_u \colon f\in\Aut({\goth M}),\; u \in C_m \right\}.
  \end{equation}
  Consequently, 
    $\Aut({\PolPap(m,{\goth M})}) \cong \Aut({\goth M})\oplus C_m$. 
\end{thm}

The case $m = 3$ is somehow exceptional in studying structures of the form
$\PolPap(m,{\goth M})$. 
Note, first, that \ref{lem:pom} and, after that, \ref{thm:autpol} do not remain valid
for $m = 3$.
An elementary reasoning shows the following
\begin{rem}
  Let $S$ be the point set of $\goth T$ 
  and $F$ be a bijection of $S\times C_3$ 
  such that for each $i\in C_3$
  and each point $x$ of $\goth T$ there is $j$ with $F(x,i) = (x,j)$.
  Then $F\in\Aut({\PolPap(3,{\goth T})})$  iff 
  either $F = \id\times\tau_u$ for some $u\in C_3$ or
  $F$ is defined by one of the following formulas%
 \footnote{%
   \begin{tabular}{r||c} $F$ & $j'$ \\ \hline\hline $u$ & $j''$ \end{tabular}
   is read as $F(u,j') = (u,j'')$}%
  :
  \begin{ctext}
    \begin{tabular}{r||c|c|c}
      $F$ & $i$ & $i+1$ & $i + 2$ 
      \\ \hline\hline
      $x$ & $j$ & $\tau_2(j)$ & $\tau_2\tau_2(j)$
      \\ \hline
      $y$ & $j+1$ & $\tau_2(j+1)$ & $\tau_2\tau_2(j+1)$
      \\ \hline
      $z$ & $j+1$ & $\tau_2(j+1)$ & $\tau_2\tau_2(j+1)$
      \\ \hline
    \end{tabular}
  \quad
    \begin{tabular}{r||c|c|c}
      $F$ & $i$ & $i+1$ & $i + 2$ 
      \\ \hline\hline
      $x$ & $j$ & $\tau_1(j)$ & $\tau_1\tau_1(j)$
      \\ \hline
      $y$ & $j+1$ & $\tau_1(j+1)$ & $\tau_1\tau_1(j+1)$
      \\ \hline
      $z$ & $j+2$ & $\tau_1(j+2)$ & $\tau_2\tau_1(j+2)$
      \\ \hline
    \end{tabular}
  \end{ctext}
  for some $i,j \in C_3$ and some labeling $x,y,z$ of the points of $\goth T$.
\end{rem}
From this we get
\begin{exm}
  Let $\goth V$ be the Veblen configuration with the points 
  $\{ a,b,c,d,p,q \}$; assume that the points $p,q$ are noncollinear in $\goth V$.
  The bijection $F$ of the points of $\PolPap(3,{\goth V})$ defined by the 
  formula 
  \begin{ctext}
    \begin{tabular}{r||c|c|c}
      F   & $0$ & $1$ & $2$ \\ \hline\hline
      $x$ & $0$ & $2$ & $1$ \\ \hline
      $y$ & $1$ & $0$ & $2$
    \end{tabular}
  \quad
  for $x\in \{ p,q \}$ and $y \in \{ a,b,c,d \}$
  \end{ctext}
  is an automorphism of $\PolPap(3,{\goth V})$.
  It is seen that $\alpha_F = \id$ but $F$ does not have form
  required in \ref{thm:autpol}.
\end{exm}

\section{Weaving and other product constructions: interrelations and applications}

In many cases a $3$-weaved configuration is a convolution with the $C_3$-group.
\begin{prop}\label{prop:hyper:polToconv0}
  Assume that a partial Steiner triple system $\goth M$ contains a hyperplane
  which is an anti-clique. 
  Let $\varepsilon_1,\varepsilon_2$ be any two elements of a group $\sf G$.
  Then 
  ${\goth M}\razy_{\varepsilon_1}{\sf G} \cong {\goth M}\razy_{\varepsilon_2}{\sf G}$.
\end{prop}
\begin{proof}
  Let $\cal H$ be a hyperplane of ${\goth M}=\struct{S,\lines}$ that is an anti-clique.
  This means the following:
  \begin{ctext}
    each line of $\goth M$ has exactly one point in common with $\cal H$.
  \end{ctext}
  Set $\varepsilon_0 = \varepsilon_2 - \varepsilon_1$ 
  (computed in ${\sf G} = \struct{G,0,+}$)
  and define the map 
  \begin{ctext}
    $\vartheta\colon S\times G \longrightarrow S\times G$,
    $\vartheta(a,i) = \left\{\begin{array}{ll}
    (a,i+\varepsilon_0) & \text{when }a\in{\cal H} \\
    (a,i)  & \text{when }a\notin{\cal H}
    \end{array}\right.$ 
    for $a\in S$, $i \in G$.
  \end{ctext}
  It is seen that $\vartheta$ maps the elements of $\lines_{\varepsilon_1}$ 
  onto the elements of $\lines_{\varepsilon_2}$ and thus it is an isomorphism of 
  ${\goth M}\razy_{\varepsilon_1} {\sf G}$ onto ${\goth M}\razy_{\varepsilon_2} {\sf G}$.
\end{proof}
Immediate from \eqref{eq:polASconv} and \ref{prop:hyper:polToconv0} is the following.
\begin{cor}\label{prop:hyper:polToconv}
  Let a partial Steiner triple system $\goth M$ contain a hyperplane
  which is an anti-clique. Then
  $\PolPap(3,{\goth M}) \cong {\goth M}\razy C_3$.
\end{cor}
\begin{fact}\label{fct:czyhipcia}
  Let $\cal H$ be an anti-clique of a 
  $\konftyp(\pointnumb,\rank,\linenumb,3)$-configuration $\goth M$.
  Then $\cal H$ is a hyperplane of $\goth M$ iff $\rank\cdot|{\cal H}| = \linenumb$
  (equivalently: iff $3\cdot|{\cal H}| = \pointnumb$).
\end{fact}
\begin{cor}\label{cor:polTOconv}
  Let $\goth M$ be one of the following partial Steiner triple systems:
  \begin{enumerate}[\rm(a)]
  \itemsep-2pt
  \item\label{zhipa:1}
    the Veblen configuration;
  \item\label{zhipa:2}  
    the Pappus Configuration or, more generally, an affine slit space
    (cf. \cite{convol}, \cite{havlik}) over $GF(3)$
    i.e. an affine space $AG(n,3)$
    with the lines parallel to a fixed affine hyperplane deleted;
  \item\label{zhipa:4}
    the configurations $\Pi^m_{\sigma_0}$ and $\Pi^m_{\id}$ for arbitrary
    $m \geq 3$.
  \end{enumerate}
  Then
    $\PolPap(3,{\goth M}) \cong {\goth M}\razy C_3$.
\end{cor}
\begin{proof}
  In case \eqref{zhipa:1} each pair of noncollinear points of the Veblen
  configuration is an anti-clique and a hyperplane.
  In case \eqref{zhipa:2} we let $\cal H$ be any hyperplane such that $\goth M$
  does not contain lines parallel to it. Then $\cal H$ is an anti-clique and a
  hyperplane in $\goth M$.
  In case \eqref{zhipa:4} the set ${\cal H} = \left\{ (0,i)\colon i =0,...,m-1 \right\}$
  is an $m$-element anti-clique and, by \ref{fct:czyhipcia},  it is a hyperplane as well.
  In each case we apply \ref{prop:hyper:polToconv} to get the claim.
\end{proof}
\begin{fact}
  Assume that a partial Steiner triple system $\goth M$ contains a hyperplane
  that is an anti-clique. 
  Then $\PolPap(m,{\goth M})$ also contains such a hyperplane%
 \footnote{%
  Analogous statement is valid for any convolution 
  ${\goth M}\razy_\varepsilon{\sf G}$.}%
  .
\end{fact}
\begin{proof}
  Let $\cal H$ be a respective hyperplane in $\goth M$.
  From \ref{lem:colin} we get that ${\cal H}\times C_m$ is an anti-clique and
  from \ref{fct:czyhipcia} it is a required hyperplane.
\end{proof}

It is not the case that each $3$-weaved configuration  
is a convolution, though.
Recall (cf. \cite{perspect}, \cite{maldem}) that
$\GrasSpace(X,2)$ with $|X| = 4$ is the Veblen configuration.
\begin{rem}\label{rem:grasNOTconv}
  Let $|X| > 4$. Then (cf. footnote \ref{foot:gras})
  $\PolPap(3,{\GrasSpace(X,2)}) \not\cong {\GrasSpace(X,2)}\razy C_3$.
\end{rem}
\begin{proof}
  Let $|X| > 4$.
  It is known that $\GrasSpace(X,2)$ contains a Desargues subconfiguration
  (cf. \cite{perspect}). 
  From \cite{convol}, 
  ${\GrasSpace(X,2)}\razy_0 C_3$ contains a Desargues subconfiguration, while
  (cf. \ref{prop:des}) $\PolPap(3,{\GrasSpace(X,2)})$
  does not contain any Desargues subconfiguration.
\end{proof}
\begin{rem}\normalfont\small
  Clearly, each point of $\goth T$ is a hyperplane. 
  By \ref{prop:hyper:polToconv0}, 
  ${\goth T}\razy_\varepsilon C_4 \cong {\goth T}\razy_0 C_4$
  for each $\varepsilon\in C_4$.
  The {\em anti-Reye configuration} 
  ${\goth T}\razy C_4$ (cf. \cite{convol}, \cite{hilbert})
  is a $\konftyp(12,4,16,3)$-configuration
  and $\PolPap(4,{\goth T})$ is a $\konftyp(12,3,12,3)$-configuration,
  and thus  ${\goth T}\razy C_4 \not\cong \PolPap(4,{\goth T})$ .
  They are also distinct in a bit stronger meaning:
  It is impossible to embed $\PolPap(4,{\goth T})$ into ${\goth T}\razy C_4$.
  Indeed, let $a$ be a point of $\goth T$.
  Two disjoint pairs of lines of ${\goth T}\razy C_4$ through $(a,0)$ 
  yield two Veblen configurations.
  Suppose $\PolPap(4,{\goth T})$ is embedded. Then two of its lines
  through $(a,0)$ should be a pair which yields a Veblen configuration,
  which is impossible, as $\PolPap(4,{\goth T})$ does not contain any
  Veblen configuration (cf. \ref{cor:vebleny}).
  Consequently, {\em the statement 
  ``under assumptions of \ref{prop:hyper:polToconv} the structure
  $\PolPap(m,{\goth M})$ is embeddable into ${\goth M}\razy C_m$''
  is not valid for $m > 3$}
 \footnote{%
  Analyzing possible series of cyclically inscribed triangles contained
  in ${\goth T}\razy C_m$ we can prove that
  {\em $\PolPap(m,{\goth T})$ is not embeddable into ${\goth T} \razy C_m$
  for any integer $m > 3$}.}%
 .
\end{rem}

As a by-product of \ref{cor:polTOconv} we get an embedding theorem
\begin{prop}
  Let $\goth M$ be 
  an affine slit space over $GF(3)$. 
  Then $\PolPap(3,{\goth M})$ can be embedded into the affine space
  over $GF(3)$.
\end{prop}
\begin{proof}
  From \ref{cor:polTOconv}, ${\goth B} := \PolPap(3,{\goth M})$ is isomorphic to 
  ${\goth M}\razy_0 C_3$. On the other hand, clearly,
  $\goth M$ is a substructure of an affine space $AG(n,3)$ for some integer $n$
  and thus $\goth B$ is a substructure of $AG(n,3)\razy C_3$.
  From \cite{convol}, $AG(n,3)\razy C_3$ is an affine slit space, embeddable into
  $AG(n+1,3)$, which proves our claim.
\end{proof}

Let ${\goth M} = \struct{S,\lines}$ be a partial Steiner triple system.
The relation defined for
$(a,i),(b,j)\in S\times C_m$ by the condition
\begin{ctext}
  $(a,i) \approx (b,j) \iff a = b$
\end{ctext}
is a congruence in $\PolPap(m,{\goth M})$ (cf. e.g. \cite{convol}).
Note also that 
  ${\PolPap(m,{\goth M})}\diagup\approx \;\cong\; {\goth M}$.

\begin{fact}\label{fct:makskliki}
  Let ${\goth M} = \struct{S,\lines}$ be a Steiner triple system. 
  The greatest   maximal anti-cliques in $\PolPap(m,{\goth M})$ 
  are the sets $\{a\}\times C_m$ with $a$ ranging over $S$, i.e. 
  the equivalence classes of the relation $\approx$.
  In case $m = 3$ these are the only maximal cliques, and 
  $\approx$ coincides with the binary non-collinearity relation.
\end{fact}
\begin{cor}\parsep-2pt
  Let ${\goth M} = \struct{S,\lines}$ be a Steiner triple system. 
  There exists the unique linear completion%
 \footnote{%
  i.e. a linear space with the point set $S$, extending $\goth M$, and 
  with the lines of the size of the lines in $\goth M$, 
  comp. \cite{linstp} (or \cite{twistfan}).}%
  \
  $\widetilde{\PolPap(3,{\goth M})}$ of $\PolPap(3,{\goth M})$.
  Its line set is the union of the family $\lines$ and the family 
  (a parallel class) $\{ \{ a \}\times C_3 \colon a \in S \}$.
\begin{rem}\label{rem:antypasz}\small
  If an STS $\goth M$ does not contain any Veblen subconfiguration
  (is {\em anti-Pasch}, or {\em Pasch-free}), 
  then (by \ref{cor:antypasz}) 
  $\PolPap(3,{\goth M})$ is also anti-Pasch. It is straightforward
  that $\widetilde{\PolPap(3,{\goth M})}$ is anti-Pasch as well.
  However, $\widetilde{\PolPap(3,{\goth M})}$ contains a miter-configuration
  for each $\goth M$.
\end{rem}
\end{cor}

The construction of a weaved configuration
has some connections with old known constructions of anti-Pasch
Steiner triple systems.
Namely, consider the group $C_3^n$ endowed with the family of blocks
$\{ u,v,2u+2v \}$ with $u,v$ ranging over pairs of distinct elements of $C_3^n$.
It is seen that that way we present the affine space $AG(n,3)$ simply.
Let a $3$-set $L = \{ a,b,c \}$ be ordered with $a < b < c$
and consider the triples
$(a,a,b)$, $(b,b,c)$, and $(c,c,a)$.
Equivalently, we can take the $C_3$-group and the family $\cal C$ defined at the
beginning of Section \ref{sec:defy} with $m = 3$.
The Bose construction (presented after \cite{griggs}, see \cite{antypasz1}, \cite{bose})
applied to the group $C_3^n$ and the set $L$ yields simply
$\widetilde{\PolPap(3,{AG(n,3)})}$.

Generally, 
the modification of the Bose construction given in \cite{griggs} and
applied to the group $C_3^k$ and the affine space $AG(n,3)$ 
yields a Steiner triple system of the parameters of the affine space $AG({n+k},3)$.
Actually, {\em it is} the affine space $AG({n+k},3)$ with the lines in one 
direction replaced by some other family of blocks.
The structure obtained by $k$-fold applying the operation
of the form 
${\goth M}\longmapsto\widetilde{\PolPap(3,{\goth M})}$
starting from $AG(n,3)$ is another (for $k > 1$) example of  an STS
with parameters of an affine space.
Since $AG(n,3)$ is anti-Pasch, the obtained structure 
is anti-Pasch as well.  
It is worth to point out that it is not an affine space, though.
\begin{prop}\label{prop:afNOTaf}
  Let $n>1$. The structure $\PolPap(3,{AG(n,3)})$
  is not embeddable to any affine space $AG(N,3)$.
  Consequently, $\widetilde{\PolPap(3,{AG(n,3)})}$
  is not embeddable to any affine space $AG(N,3)$.
  Even more generally, no structure obtained by $k$-fold applying the operation
  of the form 
  ${\goth M}\longmapsto\widetilde{\PolPap(3,{\goth M})}$
  starting from $AG(n,3)$ is embeddable as well.
  Its maximal affine subspaces are affine planes.
\end{prop}
\begin{proof}
  Let ${\goth B} = \PolPap(3,{AG(n,3)})$.
  Suppose that ${\goth B}$ is embedded into an affine 
  space ${\goth A} = AG(N,3)$ ($N>n$)
  and consider a
  triangle $\Delta = \{ (\theta,0), (b,0), (c,0) \}$
  of $\goth B$, 
  where $b\neq c,2c$, ${GF(3)}^n\ni b,c \neq \theta$ and $\theta$ is the 
  zero vector of ${GF(3)}^n$.
  From assumption, $\Delta$ spans a plane $\pi$ in 
  $\goth A$.
  We compute
  $\Delta' = \{ (2b,1), (2c,1), (2b+2c,1) \}$
  and
  $\Delta'' = \{ (b+c,2), (b+2c,2), (2b+c,2) \}$,
  and thus $\pi$ and ${\cal X} : = \Delta \cup \Delta' \cup \Delta''$ coincide.
  On the other hand the set $\cal X$
  is not a subspace in ${\goth B}$;
  indeed, $(c,0)\triop(2b,1) = (2c+b,0)$ and $(2c+b,0)\notin{\cal X}$
 \footnote{%
  The subspace spanned by $\Delta$ i.e. the
  smallest subset containing $\Delta$ and closed under $\triop$
  is the set $\pi_0\times C_3$, where $\pi_0$ is the plane spanned by
  the triangle $\{ \theta,a,b \}$ in $AG(n,3)$.
 \par
  Analogous remark remains true for arbitrary simpleks of 
  $AG(n,3)$.}%
  .  
  Thus $\pi$ is not a subspace of $\goth A$, which is a contradiction. 
  So, $\goth B$ is not embeddable, as required.
  \par\noindent
  It is straightforward that analogous reasoning applied to the triangle
  $(\theta,0_k)$, $(b,0_k)$, $(c,0_k)$ with 
  $0_k = \underbrace{0,...,0}_{k-\text{times}}$
  justifies the third nonembeddability statement.
  Computing the subspaces spanned in the considered structures by (all the possible)
  triangles we obtain our last claim.
\end{proof}

\bigskip
\begin{small}
\noindent
Author's address:\\[1ex]
Ma{\l}gorzata Pra{\.z}mowska, Krzysztof Pra{\.z}mowski
\\
Institute of Mathematics, 
University of Bia{\l}ystok
\\
\strut\quad ul. Akademicka 2, 
\strut\quad 15-267 Bia{\l}ystok, Poland
\\
e-mail: {\ttfamily malgpraz@math.uwb.edu.pl}, {\ttfamily krzypraz@math.uwb.edu.pl}
\end{small}

\end{document}